\documentclass[12pt]{article}
\pdfoutput=1
\usepackage{amsmath}    
\usepackage{graphicx}   
\usepackage{verbatim}   
\usepackage{color}      
\usepackage{subfigure}  
\usepackage{amssymb}

\newtheorem{theorem}{Theorem}[section]

\newtheorem{observation}[theorem]{Observation}

\newtheorem{corollary}[theorem]{Corollary}
\newtheorem{lemma}[theorem]{Lemma}

\newtheorem{property}[theorem]{Property}

\newcommand{\proof}{\noindent{\bf Proof.\ }}
\newcommand{\qed}{\hfill $\square$\medskip}

\begin{document}
	
\title{Strong geodetic problem  in networks: computational complexity and solution for Apollonian networks}
	
\author{
	Paul Manuel $^{a}$
	\and
	Sandi Klav\v zar $^{b,c,d}$
	\and
	Antony Xavier $^{e}$
	\and
	Andrew Arokiaraj $^{e}$
	\and
	Elizabeth Thomas $^{e}$
}

\date{}

\maketitle
\vspace{-0.8 cm}
\begin{center}
	$^a$ Department of Information Science, College of Computing Science and Engineering, Kuwait University, Kuwait \\
	{\tt pauldmanuel@gmail.com}\\
	\medskip
	
	$^b$ Faculty of Mathematics and Physics, University of Ljubljana, Slovenia\\
	{\tt sandi.klavzar@fmf.uni-lj.si}\\
	\medskip
	
	$^c$ Faculty of Natural Sciences and Mathematics, University of Maribor, Slovenia\\
	\medskip
	
	$^d$ Institute of Mathematics, Physics and Mechanics, Ljubljana, Slovenia\\
	\medskip
	
	$^e$ Department of Mathematics, Loyola College, Chennai, India \\
	
\end{center}

\begin{abstract}
The geodetic problem was introduced by Harary et al.~\cite{HLT93}. In order to model some social network problems, a similar problem is introduced in this paper and named the strong geodetic problem. The problem is solved for complete Apollonian networks. It is also proved that in general the strong geodetic problem is NP-complete. 
\end{abstract}

\noindent{\bf Keywords:} geodetic problem; strong geodetic problem; Apollonian networks;  computational complexity 

\medskip
\noindent{\bf AMS Subj.\ Class.: 05C12, 05C70, 68Q17}

\section{Introduction}
Harary et al.~\cite{HLT93} considered the following social network problem: A vertex represents  a member of the social network and an edge represents direct communication between two members of the social network. Communication among the members is restricted to only along shortest path (geodesic). Members who are lying along a geodesic are grouped together. Two coordinators supervise groups of members who lie on geodesics between the two coordinators. The problem is to identify minimum number of coordinators in such a way that each member of the social network lies on some geodesic between two coordinators. Then they modeled the above social network problem in terms of graphs as follows: Let $G=(V,E)$ be a connected graph with vertex set $V$ and edge set $E$. Let $g(x,y)$ be a geodesic between $x$ and $y$ and let $V(g(x, y))$ denote the set of vertices lying on $g(x,y)$. If $S\subseteq V$, then let $I(S)$ be the set of all geodesics between vertices of $S$ and let $V(I(S)) = \cup_{P \in I(S)} V(P)$. If $V(I(S)) = V$, then the set $S$ is called a {\em geodetic set} of $G$. The {\em geodetic problem} is to find a minimum geodetic set $S$ of $G$. 

In the seminal paper~\cite{HLT93} it was proved that the geodetic problem is NP-complete for general graphs. Moreover, Dourado et al.~\cite{DoPr10} have proved that it is NP-complete even for chordal graphs and bipartite weakly chordal graphs, while on the other hand it is polynomial on co-graphs and split graphs. Ekim et al.~\cite{EkEr12} further showed that the problem is polynomially solvable for proper interval graphs. The geodetic problem was also studied in product graphs~\cite{BrKl08}, block-cactus graphs~\cite{WaWa06}, and in line graphs~\cite{GoAs12}, while Chartrand et al.~\cite{ChZh00} investigated it in oriented graphs. Some new concepts were introduced combining geodetic and domination theory such as geodomination~\cite{ChHa01} and geodetic domination problem~\cite{HaVo10}. The hull problem which was introduced by Everett et al \cite{EvSe85} is similar to the geodetic problem. The relationship between hull problem and geodetic problem was explored by several authors~\cite{Farr05, HuTo09}. Steiner set is another concept which is similar to geodetic set. Hernando~\cite{HeJi05} and Tong~\cite{Tong09} probed the role of geodetic problem in hull and Steiner problems. For further results of the geodetic problem see~\cite{ChTo04, ChHa00, ChHa02, ChZn02} as well as the comprehensive survey~\cite{BrKo12}. 

In another situation of social networks, a set of coordinators need to be identified in such a way that each member of a social network will lie on a geodesic between two coordinators and one pair of coordinators will be able to supervise the members of only one geodesic of the social network. This situation is stronger than in the previous case. Following the geodetic problem set up, we model this social network problem as follows. Let $G=(V,E)$ be the graph corresponding to the social network. If $S\subseteq V$, then for each pair of vertices $x,y\in S$, $x\ne y$, let $\widetilde{g}(x,y)$ be a {\bf selected fixed} shortest $x,y$-path. Then we set  
$$\widetilde{I}(S)=\{\widetilde{g}(x, y) : x, y\in S\}\,,$$ 
and let $V(\widetilde{I}(S))=\cup_{P\in \widetilde{I}(S)} V(P)$. If $V(\widetilde{I}(S)) = V$ for some $\widetilde{I}(S)$, then the set $S$ is called a {\em strong geodetic set}. The {\em strong geodetic problem} is to find a minimum strong geodetic set $S$ of $G$. Clearly, the collection $\widetilde{I}(S)$ of geodesics consists of exactly $(|S|\times(|S|-1))/2$ elements. The cardinality of a minimum strong geodetic set is the {\em strong geodetic number} of $G$ and denoted by $sg(G)$. We will use $\Omega$ to denote a minimum strong geodetic set of $G$.

The rest of the paper is organized as follows. In the next section we determine the strong geodetic number of Apollonian networks while in Section~\ref{sec:NP-complete} we prove that the strong geodetic problem is NP-complete. 

\section{Apollonian networks}
\label{sec:apollonian}

In this section we discuss Apollonian networks and derive the strong geodetic number for complete Apollonian networks.  

\begin{figure}[ht!]
	\begin{center}
		\scalebox{0.5}{\includegraphics{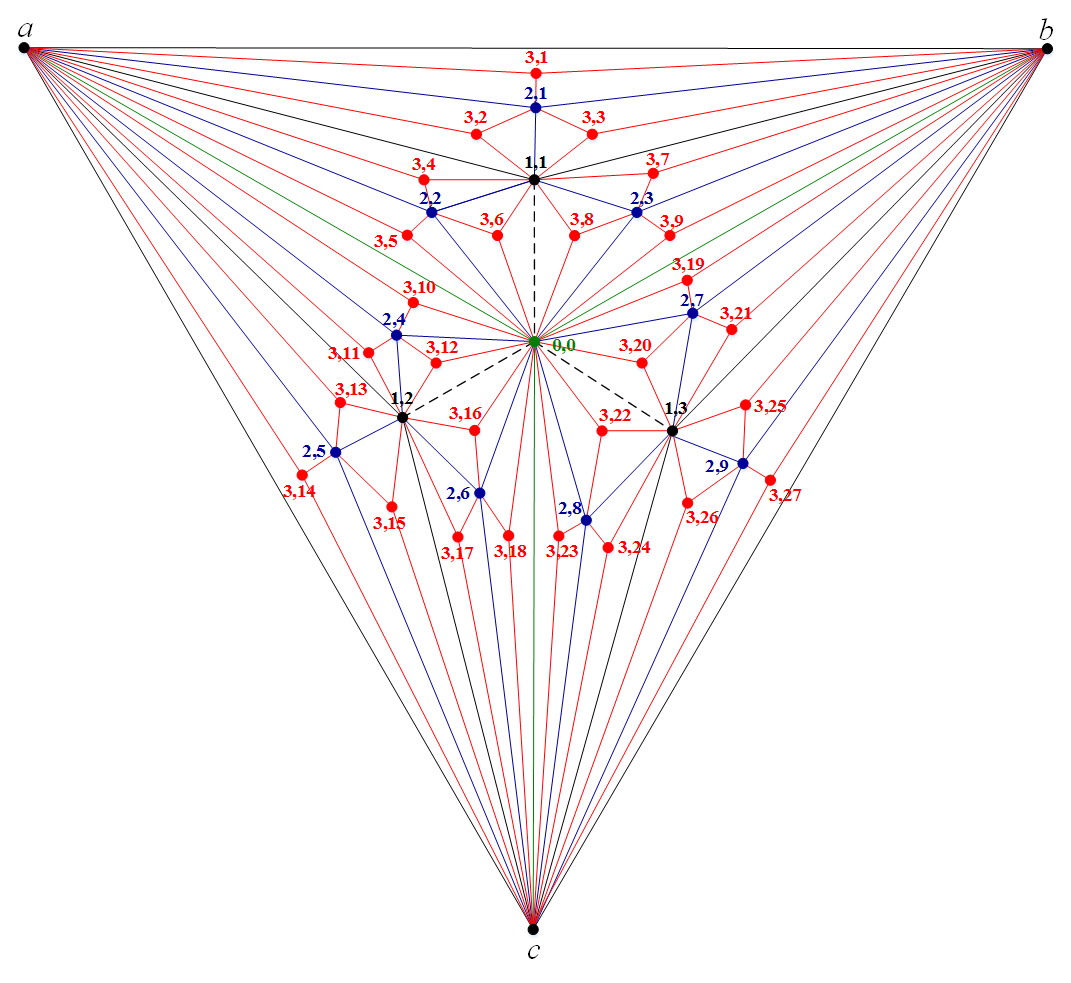}}
	\end{center}
		\caption{A complete Apollonian network $A(3)$. The 3-level vertices $T_3$ are simplicial. A vertex $(x,y)$ is written as $x,y$ due to lack of space in the diagram.}
	\label{fig:Apollonian-Networks}
\end{figure}

Apollonian networks were investigated from different points of view~\cite{GoHa75, Pell07, ZhWu14} and are constructed as follows. Start from a single triangle $\triangle(a,b,c)$. A new vertex $(0,0)$ is added inside the $\triangle(a,b,c)$ and the vertex  $(0,0)$ is connected to $a$, $b$, and $c$. The vertex $(0,0)$ is called $0$-level vertex. Inductively, $r$-level vertices are constructed from (\(r-1\))-level vertices. At $r$-th level, a new vertex $v$ is added inside a triangular face $\triangle(x,y,z)$ and the new vertex $v$ is connected to $x$, $y$, and $z$.  At each inductive step, if all (\(r-1\))-level triangular faces are filled by $r$-level vertices, then the constructed graph is called a \textit{complete Apollonian network}, otherwise it is an \textit{incomplete Apollonian network}. An $r$-level complete Apollonian network is denoted by $A(r)$. In particular, $A(0)$ is isomorphic to the complete graph $K_4$ on $4$ vertices, while $A(3)$ (together with their subgraphs $A(0)$, $A(1)$, and $A(2)$) is shown in Fig.~\ref{fig:Apollonian-Networks}. The $k$-level vertices of Apollonian network will be denoted with $(k,1),(k,2), \ldots, (k,3^k)$; in the figure we have left out the brackets in order to make the figure more transparent. 

Before stating the main result of this section we need to recall the following concepts. A vertex $v$ of a graph $G$ is {\em simplicial} if its neighborhood induces a clique. Clearly, a simplicial vertex necessarily lies in any strong geodetic set of $G$. If $G$ is a plane graph (that is, a planar graph together with a drawing in the plane), then the {\em inner dual} ${\rm inn}(G)$ of $G$ is the graph obtained by putting a vertex into each of the inner faces of $G$ and by connecting two vertices if the corresponding faces share an edge. (So the inner dual is just like the dual, except that no vertex is put into the infinite face.) Finally, we need to describe the Sierpi\'nski graphs $S_p^n$ that were  introduced in~\cite{klavzar-1997}, see the recent survey~\cite{hinz-2016+} for a wealth of information on the Sierpi\'nski graphs. More precisely, we only need the base-3 Sierpi\'nski graphs $S_3^n$ (alias Hanoi graphs $H_3^n$, see~\cite{hinz-2013}) which can be described as follows. $S_3^1 = K_3$ with $V(S_3^1) = \{0,1,2\}$. For $n\ge 2$, the Sierpi\'nski graph $S_3^{n}$ can be constructed from $3$ copies of $S_3^{n-1}$ as follows. For each $j\in \{0,1,2\}$, concatenate $j$ to the left of the vertices in a copy of $S_3^{n-1}$ and denote the obtained graph with $jS_3^{n-1}$. Then for each $i\neq j$, join copies $iS_3^{n-1}$ and $jS_3^{n-1}$ by the single edge between vertices $ij^{n-1}$ and $ji^{n-1}$. The vertices $0^n$, $1^n$, and $2^n$ are called the {\em extremal vertices} of $S_3^n$. 

In order to design an algorithm that constructs $\widetilde{I}(S)$ for Apollonian networks $A(r)$, we adopt the technique first used by Zhang, Sun, and Xu~\cite{zhang-2013} to enumerate spanning trees of Apollonian networks and then followed by Liao, Hou, and Shen~\cite{LiHu14} to calculate the Tutte polynomial. Their finding can be stated as follows. 

\begin{lemma} [\cite{LiHu14, zhang-2013}]
\label{lem:Apo-Sier} 
	If $r\ge 0$, then ${\rm inn}(A(r))$ is isomorphic to $S_3^{r+1}$.
\end{lemma}

Lemma~\ref{lem:Apo-Sier} in particular yields the following information that suits the construction of a minimum strong geodetic set of Apollonian networks. 

\begin{corollary}
\label{cor:Apo-Sier-ext} 
	Let $T_r$ denote the set of $r$-level vertices of $A(r)$. Then there is a $1$-$1$ map between the vertices of $T_r$ and the vertices of $S_3^{r}$ and, in addition, between the vertices of $V(A(r))\setminus T_r$ and the inner faces of $S_3^{r}$.
\end{corollary}

\begin{figure}[ht!]
	\begin{center}
		\scalebox{0.8}{\includegraphics{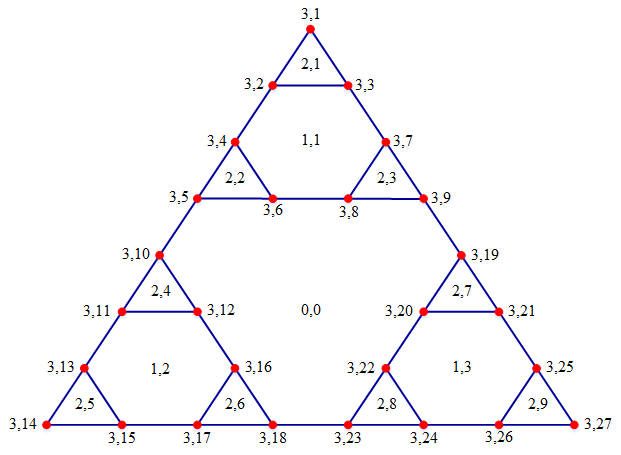}}
	\end{center}
	\caption{The Sierpi\'nski graph $S_3^{3}$. An 1-1 map between $A(3)$ and $S_3^{3}$ defined in Corollary~\ref{cor:Apo-Sier-ext} is illustrated here. The vertices of $S_3^{3}$ carry the labels of $3$-level vertices $T_3$ of $A(3)$ and the inner faces of $S_3^{3}$ carry the labels of the remaining vertices $V(A(3))\setminus T_3$ of $A(3)$. A vertex $(x,y)$ is written as $x,y$ due to lack of space in the diagram.} 
	\label{fig:Sierpinski_graph}
\end{figure}

Figs.~\ref{fig:Apollonian-Networks} and~\ref{fig:Sierpinski_graph} illustrate the $1$-$1$ map between $A(3)$ and $S_3^{3}$ which is defined in Corollary~\ref{cor:Apo-Sier-ext}.

\begin{observation}
	\label{obs:Sierpinski} 
	For each inner face $F$ of the Sierpi\'nski graph $S_3^{n}$,  there is a unique horizontal edge $xy$ of $S_3^{n}$ such that the inner face $F$ sits on the horizontal edge $xy$.
\end{observation}

For instance, in Fig.~\ref{fig:Sierpinski_graph} the inner face $(1,1)$ sits on the horizontal edge $(3,6)$-$(3,8)$ and the inner face $(2,9)$ sits on the horizontal edge $(3,26)$-$(3,27)$.

Now all is ready for the main result of this section. 

\begin{theorem} 
\label{TAPLB1} 
For Apollonian networks $A(r)$ we have $sg(A(0)) = sg(A(1)) = 4$, and if $k\ge 2$, then $sg(A(k))=3^k$. 
\end{theorem}

\proof
One can easily verify that $sg(A(0)) = sg(A(1)) = 4$. In the rest of the proof, we may thus restrict our attention to Apollonian networks $A(r)$ with $r\ge 2$. 

Let $T_r$ denote the set of $r$-level vertices of $A(r)$. Each vertex of $T_r$ is a simplicial vertex and thus it is a member of any strong geodetic set.  Thus in particular $T_{r}\subseteq \Omega$ holds, where $\Omega$ is a minimum strong geodetic set of $A(r)$. To complete the argument, it thus suffices to prove that $T_{r}$ is a strong geodetic set of $A(r)$. 

\noindent
Corollary~\ref{cor:Apo-Sier-ext} and Observation~\ref{obs:Sierpinski}  provide an algorithm to construct an $\widetilde{I}(T_r)$ of $A(r)$ as follows: 
\begin{enumerate}
	\item[(i)] 
	By Observation~\ref{obs:Sierpinski}, for each inner face $z$ of $S_3^{r}$, there is a unique horizontal edge $xy$ of $S_3^{r}$ such that the inner face $z$ sits on the edge $xy$. In other words, for each vertex $z$ of $V(A(r))\setminus T_r$ of $A(r)$, there is a unique pair of simplicial vertices $x$ and $y$ of $T_r$ such that geodesic $xzy$ of $A(r)$ covers $z$. So, we define $\widetilde{I}(T_r)$ as a collection of geodesic $xzy$ of $A(r)$ where $z$ is the inner face of $S_3^{r}$ that sits on the horizontal edge $xy$ of $S_3^{r}$.
	\item[(ii)]
	Each pair of vertices from the three extremal vertices of $S_3^{r}$ contributes a geodesic of length $2$ to $\widetilde{I}(T_r)$ to cover the vertices of $a$, $b$ and $c$. In our example of Fig.~\ref{fig:Sierpinski_graph}, (3,1), (3,14) and (3,27) are the three extremal vertices of $S_3^{r+1}$. A pair  (3,1) and (3,14) of vertices contributes a geodesic $(3,1)$-$a$-$(3,14)$ of length $2$ to $\widetilde{I}(T_r)$.  
\end{enumerate}
	
\noindent
The proof of correctness of the above algorithm is simple. As already mentioned, Corollary~\ref{cor:Apo-Sier-ext} implies that there is a 1-1 map between the vertices of $V(A(r))\setminus T_r$ and the inner faces of $S_3^{r}$. In the same way, Observation~\ref{obs:Sierpinski} implies that there is a 1-1 map between the inner faces of $S_3^{r}$ and the horizontal edges of $S_3^{r}$. Thus the above constructed $\widetilde{I}(T_r)$ consisting of geodesics of length 2 covers all the vertices of $V(A(r))\setminus T_r$.
\qed

\section{Complexity of the strong geodetic problem}
\label{sec:NP-complete}

In this section we prove that the strong geodetic problem is NP-complete. The proof's reduction will be from the dominating set problem which is a well-known NP-complete problem~\cite{GaJo90}. A set $D$ of vertices of a graph $G = (V, E)$ is a {\em dominating set} if every vertex from $V\setminus D$ has a neighbor in $D$. The {\em dominating set problem} asks whether for a given graph $G$ and integer $k$, the graph $G$ contains a dominating set of cardinality at most $k$. 

Let $G = (V,E)$ be a graph. Then construct the graph $\bar{G} = (\bar{V},\bar{E})$ as follows. The vertex set $\bar{V}$ is
\[\bar{V} = V \cup V' \cup V'',\ \textrm{where}\: V' = \{x' : x \in V\} \:\textrm{and}\: 
V'' = \{x'' : x \in V\}\,.\]
The vertex set $V'$ induces a clique and $V''$ induces an independent set. The edge set of $\bar{G}$ is $\bar{E} = E \cup E' \cup E''$, where $E'$ contains the edges of the complete graph induced by the vertices of $V'$, while $E'' = \{xx' : x \in V\} \cup  \{x'x'' : x \in V\}$.

\begin{figure}[ht!]
	\begin{center}
		\scalebox{0.4}{\includegraphics{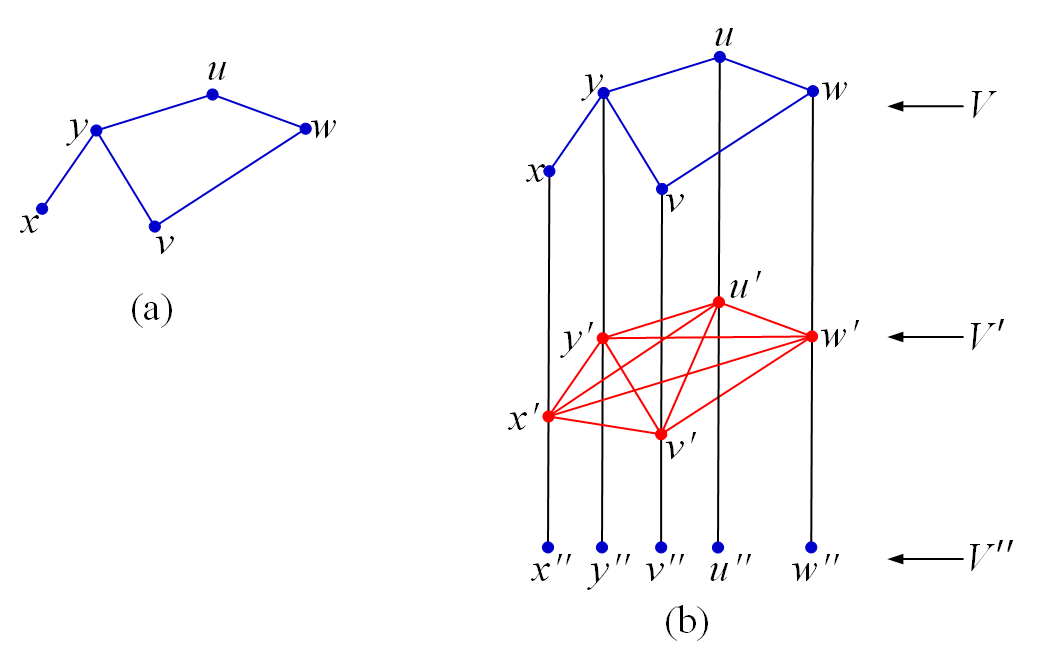}}
	\end{center}
    \caption{(a) Graph $G = (V,E)$\quad (b) $\bar{G} = (\bar{V},\bar{E})$.}
	\label{fig:NP_complete_wg}
\end{figure}

The graph $\bar{G}$ can be considered as composed of three layers: the top layer consists of $G$ itself, the middle layer forms a clique of order $|V|$, and the bottom layer is an independent set of order $|V|$. An example of the construction is presented in Fig.~\ref{fig:NP_complete_wg}.

We first observe the following fact that holds true since a pendent vertex belongs to any strong geodetic set. (Alternatively, a pendant vertex is a simplicial vertex.) 

\begin{property}
	\label{PPENWG1}
	The vertex set $V''$ of $\bar{G}$ is a subset of any strong geodetic set of $\bar{G}$. 
\end{property} 

\begin{property}
	\label{PPENWG2}
	If $X$ is a strong geodetic set of $\bar{G}$, then there exists a strong geodetic set $Y$ with $|Y| \le |X|$, such that $Y = S \cup V''$ and $S\subseteq V$. 
\end{property} 
\proof
{$X$ is a strong geodetic set of $\bar{G}$. Consider a geodesic $g(y',x)$ of $\widetilde{X}$ such that $y' \in V'$ and $x \in V$. The geodesic $g(y',x)$ is of length 2 and is of the form either $y'x'x$ or $y'yx$. The geodesic $y'x'x$ covers the vertex $x'$ of $V'$ and the geodesic $y'yx$ covers the vertex $y$ of $V$. The vertices of $V'$ are covered by geodesics $h(u'',v'')$ where $u'',v'' \in V''$ and $h(u'',v'') \in \widetilde{X}$ by Property \ref{PPENWG1}. Thus geodesic $g(y',x)$ is only of the form $y'yx$. This geodesic $y'yx$ can be replaced by $y''y'yx$ which is already in $\widetilde{X}$. Thus the vertices of $V'$ are  redundant in $X$.

Set $Y = X \setminus V'$. As discussed above, $Y$ is still a strong geodetic set of $\bar{G}$. Clearly, $Y = S \cup V''$ where $S\subseteq V$ and $|Y| \le |X|$. 
\qed

We can now prove the key fact for our reduction.   

\begin{property}
\label{PPENWG3}
$S\subseteq V$ is a dominating set of $G$ if and only if $S\cup V''$ is a strong geodetic set of $\bar{G}$.
\end{property} 

\proof
Suppose $S$ is a dominating set of $G$. Given the vertex set $S\cup V''$ in $\bar{G}$, we define the set of paths 
$$\widetilde{Y} = \{xyy'y'' : x\in S, xy\in E\}\,.$$
Note first that each path from $\widetilde{Y}$ is a geodesic. In addition, from the definition of the dominating set it easily follows that the geodesics from $\widetilde{Y}$ cover all the vertices of $V$. Next we define $$\widetilde{Z} = \{u''u'v'v'' : u'',v'' \in V''\}\,.$$
It is straightforward to observe that the geodesics from $\widetilde{Z}$ cover all the vertices of $V'\cup V''$. Now it is clear that any  
$\widetilde{I}(S\cup V'')$ that includes $\widetilde{Y} \cup \widetilde{Z}\,$
is a strong geodetic set of $\bar{G}$. 

Conversely, suppose that $S \cup V''$ is a strong geodetic set of $\bar{G}$. (We may assume that the geodetic set is of this form by Property~\ref{PPENWG2}.) Then there exists a set $\widetilde{I}(S \cup V'')$ of geodesics  such that these geodesics cover all the vertices of $\bar{G}$. Given $x \in S$ and $ y \in N(x)$, there are exactly two $x,y''$-geodesics: $xyy'y''$ and $xx'y'y''$. Among these two geodesics, only the first one $xyy'y''$ covers the vertex $y$. Therefore, $\widetilde{I}(S \cup V'')$ contains only the geodesic $xyy'y''$ among the two possible geodesics $xyy'y''$ and $xx'y'y''$.  Suppose $y\in V\setminus S$. Since $S \cup V''$ is a strong geodetic set, $y$ must be adjacent to at least one vertex from $S$, for otherwise $y$ would not be covered by $\widetilde{I}(S \cup V'')$. This implies that $S$ is a dominating set of $G$. 
\qed

Note that Property~\ref{PPENWG3} also implies that $S\subseteq V$ is a minimum dominating set of $G$ if and only if $S\cup V''$ is a minimum strong geodetic set of $\bar{G}$. Combining this result with the fact that the graph $\bar{G}$ can clearly be constructed from $G$ in a polynomial time, we have arrived at the main result of this section: 

\begin{theorem}
The strong geodetic problem is NP-complete. 
\end{theorem}

\section{Further research}
\label{sec:Fur-Research}

Even though the strong geodetic problem and the isometric path problem~\cite{Fitz99} seem to be similar, they are two different graph combinatorial problems. While the first problem minimizes the number of vertices, the second problem minimizes the number of geodesics. In this paper we have shown that the strong geodetic problem is NP-complete. To our knowledge, the complexity status of the isometric path problem is not known. Moreover, the isometric path number is known for a few graphs such as grids and block graphs but is not known even for multi-dimensional grids and other grid-like architectures. In any case, it would be useful to study the relationship between the strong geodetic problem and the isometric path problem. 



We have introduced the strong geodetic problem following the classical geodetic problem due to Harary et al.~\cite{HLT93}. We have proved that the strong geodetic problem is NP-complete. In addition, we have solved this problem for complete Apollonian networks. Further research is to investigate the strong geodetic problem for (multi-dimensional) grids, grid-like architectures, cylinders and torus. The complexity status of this problem is unknown for chordal graphs, bipartite graphs, Cayley graphs, intersection graphs, permutation graphs, etc. 

\section*{Acknowledgment}

This work was supported and funded by Kuwait University, Research Project No.\ (Q101/16).

\end{document}